\documentclass[11pt,a4paper]{article}

\usepackage[utf8]{inputenc}
\usepackage[T1]{fontenc}
\usepackage[margin=2.5cm]{geometry}
\usepackage{amsmath,amssymb,amsthm,mathtools}
\usepackage[shortlabels]{enumitem}
\usepackage{authblk}
\usepackage{xcolor}

\usepackage[colorlinks=true,linkcolor=blue,citecolor=blue,urlcolor=blue]{hyperref}

\newtheorem{theorem}{Theorem}[section]
\newtheorem{lemma}[theorem]{Lemma}
\newtheorem{proposition}[theorem]{Proposition}
\newtheorem{corollary}[theorem]{Corollary}

\theoremstyle{definition}

\newtheorem{problem}[theorem]{Problem}

\newcommand{\arxiv}[2]{\href{https://arxiv.org/abs/#1}{\texttt{arXiv:#1}} \texttt{[#2]}}
\newcommand{\doi}[1]{\url{https://doi.org/#1}}

\newcommand{\Gtriv}[1]{\Gamma^{\mathrm{triv}}_{#1}}
\newcommand{\Cay}{\operatorname{Cay}}
\newcommand{\Circ}{\operatorname{Circ}}
\newcommand{\Rcyc}{R_{\mathrm{cyc}}}
\newcommand{\Rdih}{R_{\mathrm{dih}}}
\newcommand{\Gcyc}[1]{\Gamma^{\mathrm{cyc}}_{#1}}
\newcommand{\Gdih}[1]{\Gamma^{\mathrm{dih}}_{#1}}
\newcommand{\Palt}[1]{P_{#1}^{\mathrm{alt}}}
\newcommand{\Pralt}[1]{P_{#1}^{\mathrm{ralt}}}
\newcommand{\Ssc}[1]{S_{#1}^{\mathrm{sc}}}
\newcommand{\Cmon}[1]{C_{#1}^{\mathrm{mon}}}
\newcommand{\Pmon}[1]{P_{#1}^{\mathrm{mon}}}
\newcommand{\Mnest}[1]{M_{#1}^{\mathrm{nest}}}

\title{On the cyclic and dihedral Ramsey numbers\\ of separated caterpillars and nested matchings\thanks{Email addresses:
\texttt{ivan.damnjanovic@turing.mi.sanu.ac.rs} (I.~Damnjanovi\'c),
\texttt{daniil.f.research@gmail.\allowbreak com} (D.~Furmanov),
\texttt{jon.maccarello@gmail.com} (J.~A.~Maccarello),
\texttt{julienmenet57@gmail.com} (J.~Menet).
}}

\author[1,2,3]{Ivan Damnjanović\thanks{Corresponding author.}}
\author[4]{Daniil Furmanov}
\author[5]{Jon A.~Maccarello}
\author[6]{Julien Menet}

\affil[1]{Mathematical Institute of the Serbian Academy of Sciences and Arts, Belgrade, Serbia}
\affil[2]{Faculty of Electronic Engineering, University of Niš, Niš, Serbia}
\affil[3]{FAMNIT, University of Primorska, Koper, Slovenia}
\affil[4]{Independent researcher, Vancouver, BC, Canada}
\affil[5]{Independent researcher, East Setauket, NY, United States}
\affil[6]{Independent researcher, Singapore, Singapore}

\date{}

\begin{document}

\maketitle

\begin{abstract}
We consider all graphs to have vertex sets of the form $\{0, 1, 2, \ldots, n-1\}$ for some $n \in \mathbb{N}$. For graphs $H_1, H_2, \ldots, H_k$, the \emph{cyclic Ramsey number} $\Rcyc(H_1, H_2, \ldots, H_k)$ is the smallest $n \in \mathbb{N}$ such that every $k$-edge-coloring of $K_n$ contains a monochromatic copy of $H_j$ in color~$j$ for some $j \in \{1, 2, \ldots, k\}$, with the vertices appearing in the same cyclic order as in $H_j$. The \emph{dihedral Ramsey number} $\Rdih(H_1, H_2, \ldots, H_k)$ is defined analogously, with the difference that the direction of the cyclic order is ignored.

For $n \geqslant 2$, a \emph{separated caterpillar} of order $n$ is a graph obtained by a cyclic shift of the vertex labels from a tree of order $n$ whose edges can be written as $u_1 v_1, u_2 v_2, \ldots, u_{n-1} v_{n-1}$ so that $u_1 \leqslant u_2 \leqslant \cdots \leqslant u_{n-1} < v_{n-1} \leqslant \cdots \leqslant v_2 \leqslant v_1$, while the trivial graph $K_1$ is the separated caterpillar of order one. For even $n \geqslant 2$, the \emph{nested matching} of order $n$ is the $1$-regular graph in which each vertex $v$ is adjacent to $n - 1 - v$.

We prove that if $H_1, H_2, \ldots, H_k$ are separated caterpillars, then $\Rcyc(H_1, H_2, \ldots, H_k)$ and $\Rdih(H_1, H_2, \ldots, H_k)$ are both equal to the Ramsey number of stars of the same orders, which was determined by Burr and Roberts. We also prove that $\Rcyc(F, H) = \Rdih(F, H) = 1 + (a-1)(b-1)$ for any graph $F$ of order $a \geqslant 2$ in which each vertex $v$ is adjacent to $(v + 1) \bmod a$ and any separated caterpillar $H$ of order $b \in \mathbb{N}$. Moreover, we determine the multicolor cyclic Ramsey numbers whose arguments are stars and nested matchings, provided that at most one of these nested matchings has order divisible by four. Some special cases of these results were recently obtained by Kucheriya, Petr and Yan.
\end{abstract}

\noindent \textbf{Keywords:} Ramsey numbers, cyclic Ramsey numbers, dihedral Ramsey numbers, separated caterpillars, nested matchings.

\medskip
\noindent \textbf{Mathematics Subject Classification (2020):} 05D10, 05C55, 05C35.

\section{Introduction}

We consider all graphs to be undirected, simple and finite, with vertex sets of the form $\{0, 1, 2, \ldots, n-1\}$ for some $n \in \mathbb{N}$. For a graph~$G$, we denote its vertex and edge sets by $V(G)$ and $E(G)$, respectively, and write $|G|$ for its order. We write $K_n$ for the complete graph on $n$ vertices. A \emph{$k$-edge-coloring} of a graph $G$ is a mapping from $E(G)$ to $\{1, 2, \ldots, k\}$, and for each $j \in \{1, 2, \ldots, k\}$, the \emph{color class} of color $j$ is the spanning subgraph of $G$ comprising the edges of color $j$.

Ramsey's theorem \cite{Ramsey1930} asserts that for any graphs $H_1, H_2, \ldots, H_k$, every $k$-edge-coloring of $K_n$ contains a monochromatic copy of $H_j$ in color $j$, for some $j \in \{1, 2, \ldots, k\}$, provided that $n \in \mathbb{N}$ is sufficiently large. The minimum such $n$ is the \emph{Ramsey number} $R(H_1, H_2, \ldots, H_k)$. If the embedding of $H_j$ is additionally required to be increasing, one obtains the \emph{ordered Ramsey number} $R_{\mathrm{ord}}(H_1, H_2, \ldots, H_k)$, whose systematic study was initiated independently by Balko, Cibulka, Král and Kynčl \cite{BalCiKralKyncl2015, BalCiKralKyncl2020}, and Conlon, Fox, Lee and Sudakov \cite{ConFoxLeeSu2017}; see also the survey \cite{Balko2026} by Balko. The \emph{cyclic Ramsey numbers} $\Rcyc(H_1, H_2, \ldots, H_k)$ were recently considered in \cite{BaDamSteSto2026}, and they arise naturally when the total ordering of the vertices is replaced with a cyclic ordering. More precisely, $\Rcyc(H_1, H_2, \ldots, H_k)$ is the smallest $n \in \mathbb{N}$ for which every $k$-edge-coloring of $K_n$ contains $H_j$ as a monochromatic subgraph of color $j$ for some $j \in \{1, 2, \ldots, k\}$, with the embedding $\varphi \colon V(H_j) \to V(K_n)$ required to be increasing up to a cyclic shift, i.e., there must exist $t \in V(H_j)$ such that $(\varphi(t), \varphi(t+1), \ldots, \varphi(|H_j|-1), \varphi(0), \varphi(1), \ldots, \varphi(t-1))$ is an increasing sequence.

In \cite{BaDamSteSto2026}, the standard, ordered and cyclic Ramsey numbers were unified within a group-theoretic framework through the introduction of permutational Ramsey numbers. Following \cite{DamDjor2026}, for graphs $G$ and $H$, and a permutation group $\Gamma$ on $V(H)$, we say that $H$ is \emph{$\Gamma$-embeddable} in $G$ if there exists a graph homomorphism from $H$ to $G$ of the form $\psi \circ \varphi$, where $\varphi \in \Gamma$ and $\psi \colon V(H) \to V(G)$ is an increasing injection. For graphs $H_1, H_2, \ldots, H_k$ and permutation groups $\Gamma_1, \Gamma_2, \ldots, \Gamma_k$ on their respective vertex sets, the \emph{permutational Ramsey number} $R(H_1^{\Gamma_1}, H_2^{\Gamma_2}, \ldots, H_k^{\Gamma_k})$ is the minimum $n \in \mathbb{N}$ such that for every $k$-edge-coloring of $K_n$ with color classes $G_1, G_2, \ldots, G_k$, there exists some $j \in \{1, 2, \ldots, k\}$ for which $H_j$ is $\Gamma_j$-embeddable in $G_j$.

In particular, the standard Ramsey numbers arise by taking each $\Gamma_j$ to be the symmetric group on $V(H_j)$, and the ordered Ramsey numbers arise by taking each $\Gamma_j$ to be the trivial group. Moreover, the cyclic Ramsey numbers arise by taking each $\Gamma_j$ to be the cyclic group generated by the cyclic shift permutation with respect to the natural vertex order, while taking each $\Gamma_j$ to be the group generated by the reflection with respect to the natural vertex order yields the reflective Ramsey numbers, which were introduced in \cite{BaDamSteSto2026}. If each $\Gamma_j$ is the dihedral group generated by both the cyclic shift permutation and the reflection, one obtains the \emph{dihedral Ramsey numbers} $\Rdih(H_1, H_2, \ldots, H_k)$, which were also introduced in \cite{BaDamSteSto2026}. In other words, the dihedral Ramsey numbers are defined in the same way as the cyclic Ramsey numbers, except that each monochromatic copy of $H_j$ may preserve or reverse its cyclic order.

A graph whose vertices are cyclically ordered is also known as a convex geometric graph. Thus, the Ramsey numbers of convex geometric graphs are the cyclic Ramsey numbers if the monochromatic copies are required to preserve the cyclic order, and the dihedral Ramsey numbers if they are also allowed to reverse it. Ramsey numbers of this type were studied before \cite{BaDamSteSto2026}; for example, Jayawardene, Rousseau and Harborth \cite{JayRouHar2017} determined the cyclic Ramsey numbers of certain paths and stars. Moreover, the Turán-type extremal theory of convex geometric graphs was developed by Brass, Károlyi and Valtr \cite{BraKaVal2003}, and we refer the reader to \cite{FuJiangKosMuVer2020, FuKosMuVer2020} for some more recent related results.

Small ordered and cyclic Ramsey numbers of various pairs of graph classes were computed in~\cite{BaDamSteSto2026}, while small reflective and dihedral Ramsey numbers for the same pairs of graph classes were computed in \cite{DamDjor2026}. Based on these computations, numerous conjectures were posed in both papers on these four types of Ramsey numbers in the two-color case, with the arguments belonging to specific graph families. These conjectures motivated the present work. We establish general results that include several of the conjectured formulas as special cases.

For $n \geqslant 2$, a graph of order $n$ is a \emph{separated caterpillar} if it can be obtained by a cyclic shift of the vertex labels from a tree whose edges can be written as $u_1 v_1, u_2 v_2, \ldots, u_{n-1} v_{n-1}$ so that $u_1 \leqslant u_2 \leqslant \cdots \leqslant u_{n-1} < v_{n-1} \leqslant \cdots \leqslant v_2 \leqslant v_1$, while $K_1$ is the separated caterpillar of order one. The same class of trees appeared earlier, in a different setting, in the work of Bousquet, Gledel, Narboni and Pierron \cite{BousGleNarPier2023} under the name nice caterpillars, while the name separated caterpillars was introduced by Bjerkevik, Kleist, Ueckerdt and Vogtenhuber \cite{BjerKleistUeVog2025}. We let $S_n$ denote any star graph of order $n$. Every such star is a separated caterpillar, and so are the alternating paths and reverse alternating paths considered in \cite{BaDamSteSto2026, DamDjor2026}; see Section~\ref{sc:preliminaries}. Moreover, for $n \geqslant 3$, the \emph{monotone cycle} $\Cmon{n}$ is the cycle graph of order $n$ with the underlying cycle $(0, 1, 2, \ldots, n-1)$, and we also let $\Cmon{2} \coloneqq K_2$. Finally, for any even $n \geqslant 2$, the \emph{nested matching} $\Mnest{n}$ is the $1$-regular graph of order $n$ in which each vertex $v$ is adjacent only to $n - 1 - v$.

Our main results are closely related to two classical results. Burr and Roberts \cite{BurrRo1973} determined all the Ramsey numbers whose arguments are stars, while Chvátal \cite{Chvatal1977} determined the Ramsey numbers of complete graphs versus trees.

\begin{theorem}[{\cite{BurrRo1973}}]\label{thm:burr_roberts}
Let $k \in \mathbb{N}$, let $n_1, n_2, \ldots, n_k \geqslant 2$, and let $t$ be the number of odd integers among $n_1, n_2, \ldots, n_k$. Then
\[
    R(S_{n_1}, S_{n_2}, \ldots, S_{n_k}) = \sum_{j=1}^k (n_j - 2) + \begin{cases}
        1, & \text{if } t \text{ is even and positive},\\
        2, & \text{otherwise}.
    \end{cases}
\]
\end{theorem}

\begin{theorem}[{\cite{Chvatal1977}}]\label{thm:chvatal}
For any $a, b \in \mathbb{N}$ and any tree $T$ of order $b$, we have
\[
    R(K_a, T) = 1 + (a-1)(b-1) .
\]
\end{theorem}

\noindent
Our main results can now be stated as follows.

\begin{theorem}\label{thm:main}
Let $k \in \mathbb{N}$, let $n_1, n_2, \ldots, n_k \geqslant 2$, and let $t$ be the number of odd integers among $n_1, n_2, \ldots, n_k$. Suppose that $H_j$ is a separated caterpillar of order $n_j$ for each $j \in \{1, 2, \ldots, k\}$. Then
\[
    \Rcyc(H_1, H_2, \ldots, H_k) = \Rdih(H_1, H_2, \ldots, H_k) = \sum_{j=1}^k (n_j - 2) + \begin{cases}
        1, & \text{if } t \text{ is even and positive},\\
        2, & \text{otherwise}.
    \end{cases}
\]
\end{theorem}

\begin{theorem}\label{thm:complete}
Let $a \geqslant 2$ and $b \in \mathbb{N}$, let $F$ be a graph of order $a$ such that $E(\Cmon{a}) \subseteq E(F)$, and let $H$ be a separated caterpillar of order $b$. Then
\[
    \Rcyc(F, H) = \Rdih(F, H) = 1 + (a-1)(b-1).
\]
\end{theorem}

\begin{theorem}\label{thm:star_nested}
Let $r, \ell \geqslant 0$ be such that $r + \ell \geqslant 1$, let $a_1, a_2, \ldots, a_r \geqslant 2$, and let $b_1, b_2, \ldots, b_\ell \geqslant 2$ be even integers such that at most one of them is divisible by four. Moreover, let $t_1$ be the number of odd integers among $a_1, a_2, \ldots, a_r$ and let $t_2$ be the number of integers among $b_1, b_2, \ldots, b_\ell$ that are divisible by four. Then
\begin{multline*}
    \Rcyc(S_{a_1}, \ldots, S_{a_r}, \Mnest{b_1}, \ldots, \Mnest{b_\ell}) = \Rdih(S_{a_1}, \ldots, S_{a_r}, \Mnest{b_1}, \ldots, \Mnest{b_\ell})\\
    = \sum_{j=1}^r (a_j - 2) + \sum_{j=1}^\ell (b_j - 2) + \begin{cases}
        2, & \text{if } t_1 = 0 \text{ and } t_2 \in \{0, 1\},\\
        1, & \text{if } t_1 \text{ is even and positive, and } t_2 \in \{0, 1\},\\
        2, & \text{if } t_1 \text{ is odd and } t_2 = 0,\\
        1, & \text{if } t_1 \text{ is odd and } t_2 = 1.
    \end{cases}
\end{multline*}
\end{theorem}

Theorem~\ref{thm:main} shows that every star in Theorem~\ref{thm:burr_roberts} can be independently replaced with an arbitrary separated caterpillar of the same order without changing the cyclic or the dihedral Ramsey number. Similarly, Theorem~\ref{thm:complete} shows that Theorem~\ref{thm:chvatal} remains valid for cyclic and dihedral Ramsey numbers whenever the tree is a separated caterpillar, even if the complete graph is replaced with any graph that contains the monotone cycle of the same order.

Concurrently with our work, Kucheriya, Petr and Yan \cite{KuPetrYan2026} confirmed nine conjectures from~\cite{BaDamSteSto2026} and five conjectures from~\cite{DamDjor2026}, including the conjectures implied by Theorems~\ref{thm:main}--\ref{thm:star_nested}. Their results on these implied conjectures concern two-color cyclic and dihedral Ramsey numbers and they are special cases of Theorems~\ref{thm:main}--\ref{thm:star_nested}. Moreover, Kucheriya et al.\ completed the determination of the cyclic Ramsey numbers $\Rcyc(\Mnest{a}, \Mnest{b})$ for all even $a, b \geqslant 2$. The present work was carried out independently of \cite{KuPetrYan2026}, which appeared while our manuscript was in preparation. The methods of the two works partially overlap, and we point out the overlaps in the corresponding sections.

The remainder of the paper is organized as follows. In Section~\ref{sc:preliminaries}, we introduce the necessary notation and recall some known results on which our proofs rely. In Section~\ref{sc:caterpillars}, we focus on separated caterpillars and prove Theorems~\ref{thm:main} and \ref{thm:complete}. Section~\ref{sc:nested} is devoted to nested matchings and contains the proof of Theorem~\ref{thm:star_nested}. Finally, Section~\ref{sc:conclusion} concludes the paper with a brief summary and directions for future research.

\section{Preliminaries}\label{sc:preliminaries}

In this section, we introduce the notation used throughout the paper, collect some basic facts on the embeddings used in the definitions of cyclic and dihedral Ramsey numbers, and describe a convenient labeling of the edges of separated caterpillars. We also recall several known results on which our proofs rely.

For $n \in \mathbb{N}$, let $\rho_n$ and $\tau_n$ be the permutations of $\{0, 1, 2, \ldots, n-1\}$ given by $\rho_n(x) = (x + 1) \bmod n$ and $\tau_n(x) = n - 1 - x$, i.e., the cyclic shift permutation and the reflection with respect to the natural vertex order. We write $\Gcyc{n} \coloneqq \langle \rho_n \rangle$ and $\Gdih{n} \coloneqq \langle \rho_n, \tau_n \rangle$, so that
\begin{align*}
    \Rcyc(H_1, H_2, \ldots, H_k) &= R\big(H_1^{\Gcyc{|H_1|}}, H_2^{\Gcyc{|H_2|}}, \ldots, H_k^{\Gcyc{|H_k|}}\big),\\
    \Rdih(H_1, H_2, \ldots, H_k) &= R\big(H_1^{\Gdih{|H_1|}}, H_2^{\Gdih{|H_2|}}, \ldots, H_k^{\Gdih{|H_k|}}\big).
\end{align*}
The following inequality is immediate from the definitions of cyclic and dihedral Ramsey numbers.

\begin{proposition}\label{prop:dih_cyc}
    For any $k \in \mathbb{N}$ and graphs $H_1, H_2, \ldots, H_k$, we have
    \[
        \Rdih(H_1, H_2, \ldots, H_k) \leqslant \Rcyc(H_1, H_2, \ldots, H_k) .
    \]
\end{proposition}

For a graph $H$ and a permutation $\pi$ of $V(H)$, we let $\pi(H)$ denote the graph of order $|H|$ with the edge set $\big\{ \{\pi(u), \pi(v)\} : \{u, v\} \in E(H) \big\}$. As in \cite{BaDamSteSto2026}, two graphs $H_1$ and $H_2$ of order~$n$ are \emph{rotationally isomorphic} if $H_2 = \rho_n^t(H_1)$ holds for some $t \in \mathbb{Z}$. We will frequently use the following simple facts, which are immediate from the definitions.

\begin{lemma}\label{lemma:embeddability}
Let $G$ and $H$ be graphs with $n = |H|$.
\begin{enumerate}[label=\textbf{(\roman*)}]
    \item The graph $H$ is $\Gcyc{n}$-embeddable in $G$ if and only if there exists an injective homomorphism from $H$ to $G$ that is increasing up to a cyclic shift.
    
    \item If $F$ is rotationally isomorphic to $H$, then $F$ is $\Gcyc{n}$-embeddable in $G$ if and only if $H$ is $\Gcyc{n}$-embeddable in $G$.

    \item The graph $H$ is $\Gdih{n}$-embeddable in $G$ if and only if $H$ is $\Gcyc{n}$-embeddable in $G$ or $\tau_n(H)$ is $\Gcyc{n}$-embeddable in $G$.
\end{enumerate}
\end{lemma}

For $n \in \mathbb{N}$ and $u, v \in \{0, 1, 2, \ldots, n-1\}$, we write
\[
    \delta_n(u, v) \coloneqq \min\{(u - v) \bmod n, (v - u) \bmod n\}
\]
for the circular distance between $u$ and $v$. We will also use the following elementary description of injections that are increasing up to a cyclic shift, whose straightforward proof we omit.

\begin{lemma}\label{lemma:cyclic_increasing}
Let $a, n \in \mathbb{N}$ and let $\varphi \colon \{0, 1, 2, \ldots, a-1\} \to \{0, 1, 2, \ldots, n-1\}$ be an injection. Then $\varphi$ is increasing up to a cyclic shift if and only if the map $x \mapsto (\varphi(x) - \varphi(0)) \bmod n$ is increasing. Moreover, in this case, we have
\[
    (\varphi(v) - \varphi(u)) \bmod n = (\varphi(v) - \varphi(0)) \bmod n - (\varphi(u) - \varphi(0)) \bmod n
\]
for any $0 \leqslant u < v \leqslant a-1$.
\end{lemma}

We now turn to separated caterpillars. By Lemma~\ref{lemma:embeddability}, replacing a graph of order $n$ with a rotationally isomorphic graph does not change its $\Gcyc{n}$-embeddability or $\Gdih{n}$-embeddability in any graph. In particular, the cyclic and dihedral Ramsey numbers do not change if their arguments are replaced with rotationally isomorphic graphs. Since every separated caterpillar is rotationally isomorphic to a tree of the form described in its definition, we assume without loss of generality that every separated caterpillar $H$ of order $a \geqslant 2$ is a tree with the edges $u_1 v_1, u_2 v_2, \ldots, u_{a-1} v_{a-1}$, where
\begin{equation}\label{eq:separated}
    u_1 \leqslant u_2 \leqslant \cdots \leqslant u_{a-1} < v_{a-1} \leqslant \cdots \leqslant v_2 \leqslant v_1.
\end{equation}
In this case, since the edges are distinct, $v_1 - u_1 > v_2 - u_2 > \cdots > v_{a-1} - u_{a-1}$, and since all these differences belong to $\{1, 2, \ldots, a-1\}$, we obtain $v_j - u_j = a - j$ for each $j \in \{1, 2, \ldots, a-1\}$. In particular, we have $u_1 = 0$ and $v_1 = a - 1$. Moreover, for each $j \in \{1, 2, \ldots, a-2\}$, we have either $(u_{j+1}, v_{j+1}) = (u_j + 1, v_j)$ or $(u_{j+1}, v_{j+1}) = (u_j, v_j - 1)$. Thus, the edges $u_j v_j$ and $u_{j+1} v_{j+1}$ have a common endpoint.

It is straightforward to check that if $H$ is a separated caterpillar of order $n$, then so is $\tau_n(H)$. Together with Lemma~\ref{lemma:embeddability}(iii), this gives the following lemma.

\begin{lemma}\label{lemma:cyc_to_dih}
Let $n \in \mathbb{N}$ and let $G$ be a graph such that no separated caterpillar of order $n$ is $\Gcyc{n}$-embeddable in $G$. Then no separated caterpillar of order $n$ is $\Gdih{n}$-embeddable in $G$.
\end{lemma}

For any $n \in \mathbb{N}$, the \emph{alternating path} $\Palt{n}$ is the path graph of order $n$ with the underlying path $(0, n-1, 1, n-2, 2, n-3, \ldots, \lfloor \frac{n}{2} \rfloor)$. The ordered Ramsey numbers of alternating paths were studied by Balko, Cibulka, Král and Kynčl \cite{BalCiKralKyncl2015, BalCiKralKyncl2020}, by Poljak \cite{Poljak2020}, and more recently by Kucheriya, Lo, Petr, Sgueglia and Yan \cite{KuLoPetrSgueYan2026}. Following \cite{BaDamSteSto2026}, the \emph{reverse alternating path} is $\Pralt{n} \coloneqq \tau_n(\Palt{n})$; see also \cite{KuLoPetrSgueYan2026}. Moreover, as in \cite{BaDamSteSto2026, DamDjor2026}, the \emph{start-central star} $\Ssc{n}$ is the star graph of order $n$ in which vertex $0$ is adjacent to all the other vertices, and the \emph{monotone path} $P_n^{\mathrm{mon}}$ is the path graph of order $n$ with the underlying path $(0, 1, 2, \ldots, n-1)$. For $n \geqslant 2$, the edges of $\Ssc{n}$ are $u_j v_j$ with $u_j \coloneqq 0$ and $v_j \coloneqq n - j$ for $j \in \{1, 2, \ldots, n-1\}$, while the edges of $\Palt{n}$, listed along its underlying path, are $u_j v_j$ with $u_j \coloneqq \lfloor \frac{j}{2} \rfloor$ and $v_j \coloneqq n - \lceil \frac{j}{2} \rceil$ for $j \in \{1, 2, \ldots, n-1\}$. In both cases, \eqref{eq:separated} holds with $n$ in place of $a$. Hence, $\Ssc{n}$ and $\Palt{n}$ are separated caterpillars of order~$n$, and so are $\Pralt{n}$ and every star graph of order $n$, since every star is rotationally isomorphic to a start-central star.

For any $n \in \mathbb{N}$, let $\Gtriv{n}$ denote the trivial permutation group on $\{0, 1, 2, \ldots, n-1\}$. We will derive the lower bound in Theorem~\ref{thm:complete} from the following result of \cite{DamDjor2026} on permutational Ramsey numbers in which one of the arguments is a monotone path equipped with the trivial permutation group.

\begin{theorem}[{\cite[Theorem~4.1]{DamDjor2026}}]\label{thm:mon_path}
Let $a, b \in \mathbb{N}$, let $H$ be a connected graph of order $b$, and let $\Gamma$ be a permutation group on $V(H)$. Then
\[
    R\big((P_a^{\mathrm{mon}})^{\Gtriv{a}}, H^{\Gamma}\big) = 1 + (a-1)(b-1).
\]
\end{theorem}

The upper bounds in Theorems~\ref{thm:main} and \ref{thm:complete} rely on Turán-type results. For any $n, k \in \mathbb{N}$, the Turán graph $T(n, k)$ is the complete $k$-partite graph of order $n$ in which each part has $\lfloor \frac{n}{k} \rfloor$ or $\lceil \frac{n}{k} \rceil$ vertices. We will use the classical theorem of Turán \cite{Turan1941}.

\begin{theorem}[{\cite{Turan1941}}]\label{thm:turan_classical}
Let $n, k \in \mathbb{N}$. Then, up to isomorphism, $T(n, k)$ is the unique graph with the maximum number of edges among all the graphs of order $n$ that do not contain $K_{k+1}$ as a subgraph.
\end{theorem}

We will also use a Turán-type result stating that every separated caterpillar of order $a \geqslant 2$ is $\Gcyc{a}$-embeddable in every graph of order $n$ with more than $\frac{(a-2)n}{2}$ edges. This result is a special case of a theorem of Brass, Károlyi and Valtr \cite[Theorem~7]{BraKaVal2003}, which extends an earlier result of Perles \cite[p.~292]{MoserPach1993}. Indeed, up to rotational isomorphism, the separated caterpillars are precisely the trees covered by \cite[Theorem~7]{BraKaVal2003}. Kucheriya et al.~\cite[Lemma~2.1]{KuPetrYan2026} also used this result to obtain the upper bounds in their work, and we use it in the same way in Section~\ref{sc:caterpillars}.

\begin{theorem}[{\cite[Theorem~7]{BraKaVal2003}}]\label{thm:turan}
Let $a \geqslant 2$ and $n \in \mathbb{N}$, let $H$ be a separated caterpillar of order $a$, and let $G$ be a graph of order $n$ such that $H$ is not $\Gcyc{a}$-embeddable in $G$. Then $|E(G)| \leqslant \lfloor \frac{(a-2)n}{2} \rfloor$.
\end{theorem}

\noindent
By \cite[Theorem~7]{BraKaVal2003}, the bound in Theorem~\ref{thm:turan} is the best possible whenever $n \geqslant a - 1$.

Finally, we recall some facts on circulant graphs, from which our edge colorings are built. For any $n \in \mathbb{N}$ and any $S \subseteq \{1, 2, \ldots, \lfloor \frac{n}{2} \rfloor\}$, the \emph{circulant graph} $\Circ(n, S)$ is the graph of order~$n$ in which any two vertices $u$ and $v$ are adjacent if and only if $\delta_n(u, v) \in S$. Clearly, $\Circ(n, S)$ is a regular graph whose degree equals $2|S| - 1$ if $n$ is even and $\frac{n}{2} \in S$, and $2|S|$ otherwise. In particular, for even $n$, the edge set of $\Circ(n, \{\frac{n}{2}\})$ is a perfect matching. Moreover, if $n$ is even and $x \in \mathbb{N}$ is odd with $x < \frac{n}{2}$, then $\Circ(n, \{x\})$ is the disjoint union of cycles of length $\frac{n}{\gcd(n, x)}$, which is even, so the edge set of $\Circ(n, \{x\})$ is the union of two disjoint perfect matchings.

A graph is \emph{$1$-factorizable} if its edge set can be partitioned into perfect matchings. For a finite group $\Gamma$ with the identity $e$, a \emph{connection set} of $\Gamma$ is any set $S \subseteq \Gamma \setminus \{e\}$ that is closed under inversion. For such a set $S$, the \emph{Cayley graph} $\Cay(\Gamma, S)$ is the graph with the vertex set $\Gamma$ and the edge set $\{ \{v, sv\} : v \in \Gamma, s \in S \}$. Stong \cite{Stong1985} proved the following result on the $1$-factorizability of Cayley graphs.

\begin{theorem}[{\cite{Stong1985}}]\label{thm:stong}
Let $\Gamma$ be an abelian group of even order and let $S$ be a connection set of $\Gamma$ such that $\Cay(\Gamma, S)$ is connected. Then $\Cay(\Gamma, S)$ is $1$-factorizable.
\end{theorem}

\noindent
We will use the following immediate corollary of Theorem~\ref{thm:stong}.

\begin{corollary}\label{cor:stong}
Let $n \in \mathbb{N}$ be even and let $S \subseteq \{1, 2, \ldots, \frac{n}{2}\}$ be such that $1 \in S$. Then $\Circ(n, S)$ is $1$-factorizable.
\end{corollary}

\section{Results on separated caterpillars}\label{sc:caterpillars}

In this section, we prove Theorems~\ref{thm:main} and \ref{thm:complete}. The upper bounds rely on Theorem~\ref{thm:turan}, while the lower bound in Theorem~\ref{thm:main} follows from an edge coloring built from circulant graphs. We start with a lemma on the embeddings of separated caterpillars in circulant graphs.

\begin{lemma}\label{lemma:spans}
Let $a, n \in \mathbb{N}$ with $a \geqslant 2$, let $H$ be a separated caterpillar of order $a$, let $\varphi \colon V(H) \to \{0, 1, 2, \ldots, n-1\}$ be an injection that is increasing up to a cyclic shift, and let $S \subseteq \{1, 2, \ldots, \lfloor \frac{n}{2} \rfloor\}$. Then the number of edges $uv$ of $H$ such that $\{\varphi(u), \varphi(v)\} \in E(\Circ(n, S))$ is at most the degree of $\Circ(n, S)$.
\end{lemma}

\begin{proof}
Let the edges of $H$ be labeled as in \eqref{eq:separated}, and let $f \colon V(H) \to \{0, 1, 2, \ldots, n-1\}$ be given by $f(x) \coloneqq (\varphi(x) - \varphi(0)) \bmod n$. For each $j \in \{1, 2, \ldots, a-1\}$, let $s_j \coloneqq f(v_j) - f(u_j)$. By Lemma~\ref{lemma:cyclic_increasing}, the map $f$ is increasing, and $s_j = (\varphi(v_j) - \varphi(u_j)) \bmod n$ holds for each $j \in \{1, 2, \ldots, a-1\}$. Thus, $(\varphi(u_j) - \varphi(v_j)) \bmod n = n - s_j$, so $\delta_n(\varphi(u_j), \varphi(v_j)) = \min\{s_j, n - s_j\}$. Recall that, for each $j \in \{1, 2, \ldots, a-2\}$, we have $u_j \leqslant u_{j+1}$ and $v_{j+1} \leqslant v_j$, with exactly one of these inequalities being strict. Hence, $s_j > s_{j+1}$, and the numbers $s_1, s_2, \ldots, s_{a-1}$ are mutually distinct. Now, let $T \coloneqq S \cup \{n - x : x \in S\}$. For each $j \in \{1, 2, \ldots, a-1\}$, we have $\{\varphi(u_j), \varphi(v_j)\} \in E(\Circ(n, S))$ if and only if $\min\{s_j, n - s_j\} \in S$, i.e., if and only if $s_j \in T$. Therefore, there are at most $|T|$ such $j$, and $|T|$ is the degree of $\Circ(n, S)$.
\end{proof}

For alternating paths, reverse alternating paths and stars, the case $|S| = 1$ of Lemma~\ref{lemma:spans} was also shown in \cite[Lemma~2.2]{KuPetrYan2026}. We will mostly use Lemma~\ref{lemma:spans} in the following form.

\begin{corollary}\label{cor:capacity}
Let $a, n \in \mathbb{N}$ with $a \geqslant 2$, let $H$ be a separated caterpillar of order $a$, and let $S \subseteq \{1, 2, \ldots, \lfloor \frac{n}{2} \rfloor\}$ be such that the degree of $\Circ(n, S)$ is less than $a - 1$. Then $H$ is not $\Gcyc{a}$-embeddable in $\Circ(n, S)$.
\end{corollary}

\begin{proof}
Suppose that $H$ is $\Gcyc{a}$-embeddable in $\Circ(n, S)$. By Lemma~\ref{lemma:embeddability}(i), there exists an injective homomorphism $\varphi$ from $H$ to $\Circ(n, S)$ that is increasing up to a cyclic shift. Then $\{\varphi(u), \varphi(v)\} \in E(\Circ(n, S))$ holds for all $a - 1$ edges $uv$ of $H$, which contradicts Lemma~\ref{lemma:spans}.
\end{proof}

For separated caterpillars of odd order $a \geqslant 3$, we also need the following variant of Corollary~\ref{cor:capacity}, in which an $(a-3)$-regular circulant graph is extended by a matching.

\begin{lemma}\label{lemma:saturation}
Let $a, n \in \mathbb{N}$ with $a \geqslant 3$ odd, let $H$ be a separated caterpillar of order $a$, and let $x \in \mathbb{N}$ be such that $x < \frac{n}{2}$. Suppose that $S \subseteq \{1, 2, \ldots, x - 1\}$ satisfies $|S| = \frac{a-3}{2}$ and that $F \subseteq E(\Circ(n, \{x\}))$ is a matching. Then $H$ is not $\Gcyc{a}$-embeddable in the graph $G$ obtained from $\Circ(n, S)$ by adding the edges of $F$.
\end{lemma}

\begin{proof}
Suppose that $H$ is $\Gcyc{a}$-embeddable in $G$. By Lemma~\ref{lemma:embeddability}(i), there exists an injective homomorphism $\varphi$ from $H$ to $G$ that is increasing up to a cyclic shift. Let the edges of $H$ be labeled as in \eqref{eq:separated}, and let $f$ and $s_j$, $j \in \{1, 2, \ldots, a-1\}$, be defined as in the proof of Lemma~\ref{lemma:spans}. As shown there, we have $s_1 > s_2 > \cdots > s_{a-1}$ and $\delta_n(\varphi(u_j), \varphi(v_j)) = \min\{s_j, n - s_j\}$ for each $j \in \{1, 2, \ldots, a-1\}$.

Since $|S| = \frac{a-3}{2}$ and all the elements of $S$ are smaller than $\frac{n}{2}$, the graph $\Circ(n, S)$ is $(a-3)$-regular, so Lemma~\ref{lemma:spans} implies that $\{\varphi(u), \varphi(v)\} \in E(\Circ(n, S))$ holds for at most $a - 3$ edges $uv$ of $H$. Consequently, $\{\varphi(u), \varphi(v)\} \in F$ holds for at least two of the $a - 1$ edges $uv$ of $H$. The graph $\Circ(n, \{x\})$ is $2$-regular, hence Lemma~\ref{lemma:spans} implies that there are exactly two such edges, namely, $u_h v_h$ and $u_\ell v_\ell$ with $1 \leqslant h < \ell \leqslant a - 1$. From $s_h > s_\ell$ and $\min\{s_h, n - s_h\} = \min\{s_\ell, n - s_\ell\} = x$, we obtain $s_h = n - x$ and $s_\ell = x$. If $\ell \geqslant h + 2$, then $x < s_{h+1} < n - x$, so $\delta_n(\varphi(u_{h+1}), \varphi(v_{h+1})) > x$. However, $\{\varphi(u_{h+1}), \varphi(v_{h+1})\}$ is an edge of $G$, while $\delta_n(u, v) \leqslant x$ holds for each edge $uv$ of $G$, since $S \subseteq \{1, 2, \ldots, x-1\}$ and $F \subseteq E(\Circ(n, \{x\}))$. This is a contradiction. Therefore, $\ell = h + 1$. The edges $u_h v_h$ and $u_{h+1} v_{h+1}$ have a common endpoint in $H$, so $\{\varphi(u_h), \varphi(v_h)\}$ and $\{\varphi(u_{h+1}), \varphi(v_{h+1})\}$ are two distinct edges of $F$ with a common endpoint. This is impossible, since $F$ is a matching.
\end{proof}

Corollary~\ref{cor:capacity} and Lemma~\ref{lemma:saturation} allow us to construct an edge coloring that yields the lower bound in Theorem~\ref{thm:main}.

\begin{theorem}\label{thm:coloring}
Let $k \in \mathbb{N}$, let $n_1, n_2, \ldots, n_k \geqslant 2$, let $t$ be the number of odd integers among $n_1, n_2, \ldots, n_k$, and let
\[
    n \coloneqq \sum_{j=1}^k (n_j - 2) + \begin{cases}
        0, & \text{if } t \text{ is even and positive},\\
        1, & \text{otherwise}.
    \end{cases}
\]
Then there exists a $k$-edge-coloring of $K_n$ with color classes $G_1, G_2, \ldots, G_k$ such that for each $j \in \{1, 2, \ldots, k\}$, no separated caterpillar of order $n_j$ is $\Gdih{n_j}$-embeddable in $G_j$.
\end{theorem}

\begin{proof}
By Lemma~\ref{lemma:cyc_to_dih}, it suffices to construct a $k$-edge-coloring of $K_n$ with color classes $G_1, G_2, \ldots, G_k$ such that for each $j \in \{1, 2, \ldots, k\}$, no separated caterpillar of order $n_j$ is $\Gcyc{n_j}$-embeddable in $G_j$. We consider three cases depending on the value of $t$.

\medskip
\noindent \textbf{Case 1:} $t = 0$.

In this case, all the numbers $n_1, n_2, \ldots, n_k$ are even and $n = \sum_{j=1}^k (n_j - 2) + 1$ is odd. Partition $\{1, 2, \ldots, \frac{n-1}{2}\}$ into sets $S_1, S_2, \ldots, S_k$ with $|S_j| = \frac{n_j - 2}{2}$ for each $j \in \{1, 2, \ldots, k\}$, and let $G_j \coloneqq \Circ(n, S_j)$. These circulant graphs are the color classes of a $k$-edge-coloring of $K_n$, and $G_j$ has degree $n_j - 2 < n_j - 1$ for each $j \in \{1, 2, \ldots, k\}$, so the result follows from Corollary~\ref{cor:capacity}.

\medskip
\noindent \textbf{Case 2:} $t$ is odd.

In this case, $n = \sum_{j=1}^k (n_j - 2) + 1$ is even. Let $O$ be the set of all $j \in \{1, 2, \ldots, k\}$ such that $n_j$ is odd. Choose $\beta \in O$ such that $n_\beta = \max\{n_j : j \in O\}$, let $A$ be the set of all $j \in O \setminus \{\beta\}$ with $n_j \geqslant 5$, and let $B \coloneqq O \setminus (A \cup \{\beta\})$, so that $n_j = 3$ for each $j \in B$. Partition $A$ into pairs, and if $|A|$ is odd, pair the remaining element of $A$ with an element of $B$. This yields pairs $P_1, P_2, \ldots, P_q$, each of which contains an element of $A$. Let $B_0$ be the set of all elements of $B$ that do not belong to any of these pairs, and note that $|B_0| = 2b$ for some $b \geqslant 0$, since $|A| + |B| = t - 1$ is even.

For each $\ell \in \{1, 2, \ldots, q\}$, let $d_\ell \coloneqq 1 + \sum_{j \in P_\ell} \frac{n_j - 3}{2}$, and note that $d_\ell \geqslant 2$. Moreover, let $x_\ell \coloneqq b + \sum_{j=1}^{\ell} d_j$ for $\ell \in \{0, 1, 2, \ldots, q\}$. Let $X_0 \coloneqq \{1, 2, \ldots, b\}$, and for $\ell \in \{1, 2, \ldots, q\}$, let
\[
    X_\ell \coloneqq \begin{cases}
        \{1, 2, \ldots, x_\ell\}, & \text{if } x_\ell \text{ is odd},\\
        \{1, 2, \ldots, x_\ell - 1\} \cup \{x_\ell + 1\}, & \text{if } x_\ell \text{ is even},
    \end{cases}
\]
so that $|X_\ell| = x_\ell$ and $y_\ell \coloneqq \max X_\ell$ is odd. For each $\ell \in \{1, 2, \ldots, q\}$, the inequality $d_\ell \geqslant 2$ implies $X_{\ell-1} \subseteq X_\ell$, so the set $D_\ell \coloneqq X_\ell \setminus X_{\ell-1}$ has $d_\ell$ elements and its maximum is $y_\ell$. Furthermore, a direct computation gives
\[
    \frac{n}{2} - 1 - x_q = \sum_{j \notin O} \frac{n_j - 2}{2} + \frac{n_\beta - 3}{2}.
\]
In particular, $x_q \leqslant \frac{n}{2} - 1$, and if $q \geqslant 1$, then $n_\beta \geqslant 5$, hence $y_q \leqslant x_q + 1 \leqslant \frac{n}{2} - 1$. Thus, $X_q \subseteq \{1, 2, \ldots, \frac{n}{2} - 1\}$, and we can partition the remaining $\frac{n}{2} - 1 - x_q$ elements of $\{1, 2, \ldots, \frac{n}{2} - 1\}$ into sets $S_j$, $j \notin O$, with $|S_j| = \frac{n_j - 2}{2}$, and a set $S_\beta$ with $|S_\beta| = \frac{n_\beta - 3}{2}$. Finally, for each $\ell \in \{1, 2, \ldots, q\}$, partition $D_\ell \setminus \{y_\ell\}$ into sets $S_j$, $j \in P_\ell$, with $|S_j| = \frac{n_j - 3}{2}$, and let $S_j \coloneqq \varnothing$ for each $j \in B_0$.

We now define the coloring. For each $j \in \{1, 2, \ldots, k\}$, color the edges of $\Circ(n, S_j)$ with color $j$, and color the edges of $\Circ(n, \{\frac{n}{2}\})$ with color $\beta$. For each $\ell \in \{1, 2, \ldots, q\}$, the edge set of $\Circ(n, \{y_\ell\})$ is the union of two disjoint perfect matchings, since $y_\ell$ is odd and $y_\ell < \frac{n}{2}$, and we let each color $j \in P_\ell$ receive one of them. If $b \geqslant 1$, then $\Circ(n, X_0)$ is $2b$-regular and $1 \in X_0$, so by Corollary~\ref{cor:stong}, its edge set is the union of $2b$ mutually disjoint perfect matchings,  and we let each color $j \in B_0$ receive one of them. Since the sets $X_0$, $\{y_1\}, \{y_2\}, \ldots, \{y_q\}$ and $S_1, S_2, \ldots, S_k$ partition $\{1, 2, \ldots, \frac{n}{2} - 1\}$, every edge of $K_n$ receives exactly one color.

It remains to verify the coloring. If $j \notin O$, then $G_j = \Circ(n, S_j)$ has degree $n_j - 2 < n_j - 1$, and if $j = \beta$, then $G_j = \Circ(n, S_\beta \cup \{\frac{n}{2}\})$ has degree $n_\beta - 2 < n_\beta - 1$. In both cases, the result follows from Corollary~\ref{cor:capacity}. If $j \in A$, then $j \in P_\ell$ for some $\ell \in \{1, 2, \ldots, q\}$, and $G_j$ is obtained from $\Circ(n, S_j)$ by adding a matching contained in $E(\Circ(n, \{y_\ell\}))$, where $|S_j| = \frac{n_j - 3}{2}$ and all the elements of $S_j$ are smaller than $y_\ell < \frac{n}{2}$, so the result follows from Lemma~\ref{lemma:saturation}. Finally, if $j \in B$, then $G_j$ is a perfect matching, while a separated caterpillar of order three has two edges with a common endpoint, so it admits no injective homomorphism to $G_j$.

\medskip
\noindent \textbf{Case 3:} $t$ is even and positive.

Choose $h \in \{1, 2, \ldots, k\}$ such that $n_h$ is odd, and let $n'_h \coloneqq n_h - 1$ and $n'_j \coloneqq n_j$ for each $j \in \{1, 2, \ldots, k\} \setminus \{h\}$. Then exactly $t - 1$ of the numbers $n'_1, n'_2, \ldots, n'_k$ are odd and $\sum_{j=1}^k (n'_j - 2) + 1 = n$, so Case~2, applied to $n'_1, n'_2, \ldots, n'_k$, yields a $k$-edge-coloring of $K_n$ with color classes $G_1, G_2, \ldots, G_k$ such that for each $j \in \{1, 2, \ldots, k\}$, no separated caterpillar of order $n'_j$ is $\Gcyc{n'_j}$-embeddable in $G_j$. It remains to observe that, since no separated caterpillar of order $n_h - 1$ is $\Gcyc{n_h - 1}$-embeddable in $G_h$, no separated caterpillar of order $n_h$ is $\Gcyc{n_h}$-embeddable in $G_h$ either.
\end{proof}

We emphasize that Theorem~\ref{thm:coloring} gives more than what is needed for the lower bound in Theorem~\ref{thm:main}. Indeed, the constructed edge coloring depends only on the orders $n_1, n_2, \ldots, n_k$, and for each $j \in \{1, 2, \ldots, k\}$, the color class $G_j$ simultaneously avoids all the separated caterpillars of order $n_j$, rather than a single given one. We are now ready to prove Theorem~\ref{thm:main}.

\begin{proof}[Proof of Theorem~\ref{thm:main}]
Let
\[
    n \coloneqq \sum_{j=1}^k (n_j - 2) + \begin{cases}
        1, & \text{if } t \text{ is even and positive},\\
        2, & \text{otherwise}.
    \end{cases}
\]
Theorem~\ref{thm:coloring} implies that $\Rdih(H_1, H_2, \ldots, H_k) \geqslant n$, so, by Proposition~\ref{prop:dih_cyc}, it remains to prove that $\Rcyc(H_1, H_2, \ldots, H_k) \leqslant n$.

By way of contradiction, suppose that there exists a $k$-edge-coloring of $K_n$ with color classes $G_1, G_2, \ldots, G_k$ such that $H_j$ is not $\Gcyc{n_j}$-embeddable in $G_j$ for any $j \in \{1, 2, \ldots, k\}$. By Theorem~\ref{thm:turan}, we have $|E(G_j)| \leqslant \lfloor \frac{(n_j - 2)n}{2} \rfloor$ for each $j \in \{1, 2, \ldots, k\}$. If $t = 0$ or $t$ is odd, then $n = \sum_{j=1}^k (n_j - 2) + 2$, so we obtain
\[
    \binom{n}{2} = \sum_{j=1}^k |E(G_j)| \leqslant \sum_{j = 1}^k \frac{(n_j - 2) n}{2} = \frac{n(n-2)}{2} < \binom{n}{2},
\]
which is a contradiction. On the other hand, if $t$ is even and positive, then $n = \sum_{j=1}^k (n_j - 2) + 1$ is odd. Thus, for every $j \in \{ 1, 2, \ldots, k \}$ such that $n_j$ is odd, we have $|E(G_j)| \leqslant \frac{(n_j - 2)n - 1}{2}$, hence
\[
    \binom{n}{2} = \sum_{j=1}^k |E(G_j)| \leqslant \sum_{j=1}^k \frac{(n_j - 2)n}{2} - \frac{t}{2} = \frac{n(n-1)}{2} - \frac{t}{2} < \binom{n}{2},
\]
which is again a contradiction. Therefore, $\Rcyc(H_1, H_2, \ldots, H_k) \leqslant n$.
\end{proof}

In particular, for $k = 2$, Theorem~\ref{thm:main} yields
\[
    \Rcyc(F_1, F_2) = \Rdih(F_1, F_2) = a + b - 2 - (ab \bmod 2)
\]
for any $a, b \geqslant 2$, $F_1 \in \{S_a, \Palt{a}, \Pralt{a}\}$ and $F_2 \in \{S_b, \Palt{b}, \Pralt{b}\}$. For several choices of $F_1$ and $F_2$, this formula was conjectured in \cite{BaDamSteSto2026, DamDjor2026}, and all these conjectures were recently confirmed by Kucheriya et al.~\cite{KuPetrYan2026}. For each of these choices, our proof of the upper bound in Theorem~\ref{thm:main} is essentially the same as the proof of \cite[Theorem~2.4]{KuPetrYan2026}.

We now turn to the proof of Theorem~\ref{thm:complete}.

\begin{proof}[Proof of Theorem~\ref{thm:complete}]
Let $n \coloneqq 1 + (a-1)(b-1)$. Since $\rho_a(\Cmon{a}) = \Cmon{a}$ and $\tau_a(\Cmon{a}) = \Cmon{a}$, from $E(\Cmon{a}) \subseteq E(F)$, $E(\Pmon{a}) \subseteq E(\Cmon{a})$ and Theorem~\ref{thm:mon_path}, we obtain
\begin{align*}
    \Rdih(F, H) = R\big(F^{\Gdih{a}}, H^{\Gdih{b}}\big) &\geqslant R\big((\Cmon{a})^{\Gdih{a}}, H^{\Gdih{b}}\big) = R\big((\Cmon{a})^{\Gtriv{a}}, H^{\Gdih{b}}\big)\\
    &\geqslant R\big((\Pmon{a})^{\Gtriv{a}}, H^{\Gdih{b}}\big) = n.
\end{align*}
Thus, by Proposition~\ref{prop:dih_cyc}, it remains to prove that $\Rcyc(F, H) \leqslant n$.

By way of contradiction, suppose that there exists a $2$-edge-coloring of $K_n$ with color classes $G_1$ and $G_2$ such that $F$ is not $\Gcyc{a}$-embeddable in $G_1$ and $H$ is not $\Gcyc{b}$-embeddable in $G_2$. In this case, we have $b \geqslant 2$. Moreover, $G_1$ clearly cannot contain $K_a$ as a subgraph, so, by Theorem~\ref{thm:turan_classical}, we have $|E(G_1)| \leqslant |E(T(n, a-1))|$. The Turán graph $T(n, a-1)$ has one part of size $b$ and $a - 2$ parts of size $b - 1$, so
\begin{equation}\label{eq:aux_1}
    |E(G_1)| \leqslant \binom{n}{2} - \binom{b}{2} - (a - 2) \binom{b - 1}{2}.
\end{equation}
On the other hand, Theorem~\ref{thm:turan} gives $|E(G_2)| \leqslant \frac{(b-2)n}{2} = (a-1)\binom{b-1}{2} + \frac{b-2}{2}$, which together with \eqref{eq:aux_1} implies
\[
    \binom{n}{2} = |E(G_1)| + |E(G_2)| \leqslant \binom{n}{2} - \binom{b}{2} + \binom{b-1}{2} + \frac{b-2}{2} = \binom{n}{2} - \frac{b}{2},
\]
yielding a contradiction. Therefore, $\Rcyc(F, H) \leqslant n$.
\end{proof}

For $F \in \{K_a, \Cmon{a}\}$ and $H = \Palt{b}$, the formula from Theorem~\ref{thm:complete} was conjectured in \cite{BaDamSteSto2026, DamDjor2026}, and these conjectures were recently confirmed by Kucheriya et al.~\cite{KuPetrYan2026}. For these choices of $F$ and $H$, our proof of the upper bound in Theorem~\ref{thm:complete} is essentially the same as the proof of \cite[Lemma~2.6]{KuPetrYan2026}. Moreover, since every star is a separated caterpillar, Theorem~\ref{thm:complete} also gives $\Rcyc(K_a, S_b) = \Rcyc(\Cmon{a}, S_b) = 1 + (a-1)(b-1)$ for all $a \geqslant 2$ and $b \in \mathbb{N}$, as already shown in \cite[Corollaries~4.19 and~4.26]{BaDamSteSto2026}.

\section{Results on nested matchings}\label{sc:nested}

In this section, we prove Theorem~\ref{thm:star_nested}. First, we note that any star graph or nested matching of order $n$ is $\Gdih{n}$-embeddable in a graph if and only if it is $\Gcyc{n}$-embeddable in it, so the corresponding cyclic and dihedral Ramsey numbers from Theorem~\ref{thm:star_nested} coincide. Moreover, for any $a \geqslant 2$, a star of order $a$ is $\Gcyc{a}$-embeddable in a graph $G$ if and only if $G$ has a vertex of degree at least $a - 1$.

For $n \in \mathbb{N}$ and $s \in \{0, 1, 2, \ldots, n-1\}$, let $L_n(s)$ denote the set of all the edges $\{x, y\}$ of $K_n$ such that $x + y \equiv s \pmod{n}$. Clearly, the sets $L_n(0), L_n(1), \ldots, L_n(n-1)$ partition $E(K_n)$. The following lemma bounds the number of edges in each of these sets for the graphs in which $\Mnest{b}$ is not $\Gcyc{b}$-embeddable.

\begin{lemma}\label{lemma:nested_reflection}
Let $n \in \mathbb{N}$, let $b \geqslant 2$ be even, let $s \in \{0, 1, 2, \ldots, n-1\}$, and let $G$ be a graph of order $n$ such that $\Mnest{b}$ is not $\Gcyc{b}$-embeddable in $G$. Then $|E(G) \cap L_n(s)| \leqslant \frac{b}{2} - 1$.
\end{lemma}

\begin{proof}
The statement trivially holds for $n = 1$, so we assume that $n \geqslant 2$. Let $G'$ be the graph of order $n$ with the edge set $L_n(s)$. It is straightforward to observe that $G'$ is rotationally isomorphic to one of the three graphs below, depending on the parities of $n$ and $s$:
\begin{enumerate}[label=\textbf{(\arabic*)}]
    \item the graph obtained from $\Mnest{n - 1}$ by adding an isolated vertex $n - 1$, if $n$ is odd;
    \item $\Mnest{n}$, if $n$ is even and $s$ is odd; and
    \item the graph of order $n$ with the edges $\{ v, n - 2 - v \}$ with $0 \leqslant v \leqslant \frac{n}{2} - 2$, if $n$ and $s$ are both even.
\end{enumerate}
In all three cases, if we had $|E(G) \cap L_n(s)| \geqslant \frac{b}{2}$, then $\Mnest{b}$ would be $\Gcyc{b}$-embeddable in $G$, yielding a contradiction.
\end{proof}

Lemma~\ref{lemma:nested_reflection} implies that every graph $G$ of order $n$ in which $\Mnest{b}$ is not $\Gcyc{b}$-embeddable satisfies $|E(G)| \leqslant \frac{(b-2)n}{2}$. This is the counterpart of Theorem~\ref{thm:turan} for nested matchings. The same argument was used by Kucheriya et~al.\ in the proof of \cite[Proposition~3.2]{KuPetrYan2026}, where it was also shown that this bound is the best possible whenever $n \geqslant b$. The following lemma plays the role of Corollary~\ref{cor:capacity} for nested matchings.

\begin{lemma}\label{lemma:layer_blocks}
Let $n \in \mathbb{N}$, let $b \geqslant 2$ be even, and let $S \subseteq \{1, 2, \ldots, \lfloor \frac{n}{2} \rfloor\}$ be a set of $\frac{b}{2} - 1$ consecutive integers. Suppose that $b \equiv 2 \pmod{4}$ or $\lfloor \frac{n}{2} \rfloor \in S$. Then $\Mnest{b}$ is not $\Gcyc{b}$-embeddable in $\Circ(n, S)$.
\end{lemma}

\begin{proof}
By way of contradiction, suppose that $\Mnest{b}$ is $\Gcyc{b}$-embeddable in $\Circ(n, S)$, and let $\varphi$ be an injective homomorphism from $\Mnest{b}$ to $\Circ(n, S)$ that is increasing up to a cyclic shift. Let $f \colon V(\Mnest{b}) \to \{0, 1, 2, \ldots, n-1\}$ be given by $f(x) \coloneqq (\varphi(x) - \varphi(0)) \bmod n$, and for each $j \in \{0, 1, \ldots, \frac{b}{2} - 1\}$, let $\lambda_j \coloneqq f(b - 1 - j) - f(j)$. By Lemma~\ref{lemma:cyclic_increasing}, the map $f$ is increasing, and for each $j \in \{0, 1, \ldots, \frac{b}{2} - 1\}$, we have $\lambda_j = (\varphi(b - 1 - j) - \varphi(j)) \bmod n$, so $\min\{\lambda_j, n - \lambda_j\} = \delta_n(\varphi(j), \varphi(b - 1 - j)) \in S$, i.e., $\lambda_j \in T \coloneqq S \cup \{n - x : x \in S\}$. Moreover, $\lambda_j - \lambda_{j+1} = \big(f(b - 1 - j) - f(b - 2 - j)\big) + \big(f(j + 1) - f(j)\big) \geqslant 2$ for each $j \in \{0, 1, \ldots, \frac{b}{2} - 2\}$.

If $\lfloor \frac{n}{2} \rfloor \in S$, then $T$ is a set of at most $b - 2$ consecutive integers, so it cannot contain all the $\frac{b}{2}$ numbers $\lambda_0, \lambda_1, \ldots, \lambda_{\frac{b}{2} - 1}$, since $\lambda_j - \lambda_{j+1} \geqslant 2$ for each $j \in \{0, 1, \ldots, \frac{b}{2} - 2\}$. Otherwise, $b \equiv 2 \pmod{4}$, and $T$ is the union of two sets of $\frac{b}{2} - 1$ consecutive integers, each of which contains at most $\frac{b - 2}{4}$ of these numbers. In both cases, we obtain a contradiction.
\end{proof}

We mention in passing that, for $S = \{1, 2, \ldots, \frac{b}{2} - 1\}$ with $b \equiv 2 \pmod{4}$, and for odd $n$ with $\lfloor \frac{n}{2} \rfloor \in S$, Lemma~\ref{lemma:layer_blocks} yields two of the constructions used by Kucheriya et~al.\ in the proof of \cite[Proposition~3.2]{KuPetrYan2026}, whose verification also relies on the same observation on the circular distances. We are now ready to construct the edge coloring that yields the lower bound in Theorem~\ref{thm:star_nested}.

\begin{lemma}\label{lemma:star_nested_coloring}
Let $r, \ell \geqslant 0$ be such that $r + \ell \geqslant 1$, let $a_1, a_2, \ldots, a_r \geqslant 2$, and let $b_1, b_2, \ldots, b_\ell \geqslant 2$ be even integers such that at most one of them is divisible by four. Moreover, let $t_1$ be the number of odd integers among $a_1, a_2, \ldots, a_r$, let $t_2$ be the number of integers among $b_1, b_2, \ldots, b_\ell$ that are divisible by four, and let
\[
    n \coloneqq \sum_{j=1}^r (a_j - 2) + \sum_{j=1}^\ell (b_j - 2) + \begin{cases}
        1, & \text{if } t_1 = 0 \text{ and } t_2 \in \{0, 1\},\\
        0, & \text{if } t_1 \text{ is even and positive, and } t_2 \in \{0, 1\},\\
        1, & \text{if } t_1 \text{ is odd and } t_2 = 0,\\
        0, & \text{if } t_1 \text{ is odd and } t_2 = 1.
    \end{cases}
\]
Then there exists an $(r + \ell)$-edge-coloring of $K_n$ with color classes $G_1, G_2, \ldots, G_{r+\ell}$ such that $S_{a_j}$ is not $\Gcyc{a_j}$-embeddable in $G_j$ for any $j \in \{1, 2, \ldots, r\}$, and $\Mnest{b_j}$ is not $\Gcyc{b_j}$-embeddable in $G_{r+j}$ for any $j \in \{1, 2, \ldots, \ell\}$.
\end{lemma}

\begin{proof}
If $\ell = 0$, then the statement follows from Theorem~\ref{thm:coloring}, so assume that $\ell \geqslant 1$. Moreover, we can assume without loss of generality that $b_j \equiv 2 \pmod{4}$ for each $j \in \{1, 2, \ldots, \ell - 1\}$. We split the argument into four cases depending on the values of $t_1$ and $t_2$.

\medskip
\noindent \textbf{Case 1:} $t_1 = 0$.

In this case, $n = \sum_{j=1}^r (a_j - 2) + \sum_{j=1}^\ell (b_j - 2) + 1$ is odd. Partition $\{1, 2, \ldots, \frac{n - 1}{2} \}$ into sets $S_1, S_2, \ldots, S_{r+\ell}$ of consecutive integers, placed in this order, such that $|S_j| = \frac{a_j - 2}{2}$ for each $j \in \{1, 2, \ldots, r\}$ and $|S_{r+j}| = \frac{b_j - 2}{2}$ for each $j \in \{1, 2, \ldots, \ell\}$, and let $G_j \coloneqq \Circ(n, S_j)$ for each $j \in \{1, 2, \ldots, r + \ell\}$. These circulant graphs are the color classes of an $(r + \ell)$-edge-coloring of $K_n$. For each $j \in \{1, 2, \ldots, r\}$, the graph $G_j$ is $(a_j - 2)$-regular, so $S_{a_j}$ is not $\Gcyc{a_j}$-embeddable in it. Moreover, for each $j \in \{1, 2, \ldots, \ell - 1\}$, we have $b_j \equiv 2 \pmod{4}$, while $b_\ell = 2$ or $\frac{n - 1}{2} \in S_{r+\ell}$. Thus, by Lemma~\ref{lemma:layer_blocks}, $\Mnest{b_j}$ is not $\Gcyc{b_j}$-embeddable in $G_{r+j}$ for any $j \in \{1, 2, \ldots, \ell\}$.

\medskip
\noindent \textbf{Case 2:} $t_1$ is even and positive.

In this case, $n = \sum_{j=1}^r (a_j - 2) + \sum_{j=1}^\ell (b_j - 2)$ is even. Let $p \coloneqq \sum_{j=1}^r \frac{a_j - 2}{2}$, and note that $p \geqslant 1$. As in Case~1, partition $\{p + 1, p + 2, \ldots, \frac{n}{2}\}$ into sets $S_{r+1}, S_{r+2}, \ldots, S_{r+\ell}$ of consecutive integers, placed in this order, such that $|S_{r+j}| = \frac{b_j - 2}{2}$ for each $j \in \{1, 2, \ldots, \ell\}$, and let $G_{r+j} \coloneqq \Circ(n, S_{r+j})$ for each $j \in \{1, 2, \ldots, \ell\}$. Since $b_\ell = 2$ or $\frac{n}{2} \in S_{r+\ell}$, Lemma~\ref{lemma:layer_blocks} implies that $\Mnest{b_j}$ is not $\Gcyc{b_j}$-embeddable in $G_{r+j}$ for any $j \in \{1, 2, \ldots, \ell\}$. It remains to distribute the edges of $\Circ(n, \{1, 2, \ldots, p\})$ among the colors $1, 2, \ldots, r$ so that every vertex of $G_j$ has degree at most $a_j - 2$ for each $j \in \{1, 2, \ldots, r\}$. If $p < \frac{n}{2}$, then $\Circ(n, \{1, 2, \ldots, p\})$ is $2p$-regular, so by Corollary~\ref{cor:stong}, its edge set is the union of $2p$ mutually disjoint perfect matchings, and we let each color $j \in \{1, 2, \ldots, r\}$ receive $a_j - 2$ of them. On the other hand, if $p = \frac{n}{2}$, then $\Circ(n, \{1, 2, \ldots, p\}) = K_n$ is $(2p - 1)$-regular, so its edge set is the union of $2p - 1$ mutually disjoint perfect matchings. Choose $h \in \{1, 2, \ldots, r\}$ such that $a_h$ is odd. Then we can let the color $h$ receive $a_h - 3$ of these perfect matchings, while each color $j \in \{1, 2, \ldots, r\} \setminus \{h\}$ receives $a_j - 2$ of them.

\medskip
\noindent \textbf{Case 3:} $t_1$ is odd and $t_2 = 0$.

In this case, $n = \sum_{j=1}^r (a_j - 2) + \sum_{j=1}^\ell (b_j - 2) + 1$ is even. Let $p \coloneqq \frac{1}{2} \big( \sum_{j=1}^r (a_j - 2) - 1 \big)$. As in Case~1, partition $\{p + 1, p + 2, \ldots, \frac{n}{2} - 1\}$ into sets $S_{r+1}, S_{r+2}, \ldots, S_{r+\ell}$ of consecutive integers such that $|S_{r+j}| = \frac{b_j - 2}{2}$ for each $j \in \{1, 2, \ldots, \ell\}$, and let $G_{r+j} \coloneqq \Circ(n, S_{r+j})$ for each $j \in \{1, 2, \ldots, \ell\}$. Since $t_2 = 0$, Lemma~\ref{lemma:layer_blocks} implies that $\Mnest{b_j}$ is not $\Gcyc{b_j}$-embeddable in $G_{r+j}$ for any $j \in \{1, 2, \ldots, \ell\}$. It remains to distribute the edges of $\Circ(n, \{1, 2, \ldots, p\} \cup \{\frac{n}{2}\})$ among the colors $1, 2, \ldots, r$ so that every vertex of $G_j$ has degree at most $a_j - 2$ for each $j \in \{1, 2, \ldots, r\}$. This circulant graph is $(2p + 1)$-regular, and its edge set is the union of $2p + 1$ mutually disjoint perfect matchings, by Corollary~\ref{cor:stong} if $p \geqslant 1$, and trivially if $p = 0$. Since $\sum_{j=1}^r (a_j - 2) = 2p + 1$, we can let each color $j \in \{1, 2, \ldots, r\}$ receive $a_j - 2$ of these perfect matchings.

\medskip
\noindent \textbf{Case 4:} $t_1$ is odd and $t_2 = 1$.

In this case, $n = \sum_{j=1}^r (a_j - 2) + \sum_{j=1}^\ell (b_j - 2)$ is odd. Let $p \coloneqq \frac{1}{2} \big( \sum_{j=1}^r (a_j - 2) - 1 \big)$. As in Case~1, partition $\{p + 1, p + 2, \ldots, \frac{n - 1}{2}\}$ into sets $S_{r+1}, S_{r+2}, \ldots, S_{r+\ell}$ of consecutive integers, placed in this order, such that $|S_{r+j}| = \frac{b_j - 2}{2}$ for each $j \in \{1, 2, \ldots, \ell\}$, and let $G_{r+j} \coloneqq \Circ(n, S_{r+j})$ for each $j \in \{1, 2, \ldots, \ell\}$. Since $\frac{n - 1}{2} \in S_{r+\ell}$, Lemma~\ref{lemma:layer_blocks} implies that $\Mnest{b_j}$ is not $\Gcyc{b_j}$-embeddable in $G_{r+j}$ for any $j \in \{1, 2, \ldots, \ell\}$. It remains to distribute the edges of $\Circ(n, \{1, 2, \ldots, p\})$ among the colors $1, 2, \ldots, r$ so that every vertex of $G_j$ has degree at most $a_j - 2$ for each $j \in \{1, 2, \ldots, r\}$. This circulant graph is $2p$-regular, so by Vizing's theorem \cite{Vizing1964}, its edge set is the union of $2p + 1$ mutually disjoint matchings. Since $\sum_{j=1}^r (a_j - 2) = 2p + 1$, we can let each color $j \in \{1, 2, \ldots, r\}$ receive $a_j - 2$ of these matchings.
\end{proof}

\noindent
We are now ready to prove Theorem~\ref{thm:star_nested}.

\begin{proof}[Proof of Theorem~\ref{thm:star_nested}]
Let
\[
    n \coloneqq \sum_{j=1}^r (a_j - 2) + \sum_{j=1}^\ell (b_j - 2) + \begin{cases}
        2, & \text{if } t_1 = 0 \text{ and } t_2 \in \{0, 1\},\\
        1, & \text{if } t_1 \text{ is even and positive, and } t_2 \in \{0, 1\},\\
        2, & \text{if } t_1 \text{ is odd and } t_2 = 0,\\
        1, & \text{if } t_1 \text{ is odd and } t_2 = 1.
    \end{cases}
\]
As noted at the beginning of this section, the cyclic and dihedral Ramsey numbers from Theorem~\ref{thm:star_nested} coincide, so it suffices to prove that $\Rcyc(S_{a_1}, \ldots, S_{a_r}, \Mnest{b_1}, \ldots, \Mnest{b_\ell}) = n$. Lemma~\ref{lemma:star_nested_coloring} implies that $\Rcyc(S_{a_1}, \ldots, S_{a_r}, \Mnest{b_1}, \ldots, \Mnest{b_\ell}) \geqslant n$, so it remains to prove that $\Rcyc(S_{a_1}, \ldots, S_{a_r}, \Mnest{b_1}, \ldots, \Mnest{b_\ell}) \leqslant n$.

By way of contradiction, suppose that there exists an $(r + \ell)$-edge-coloring of $K_n$ with color classes $G_1, G_2, \ldots, G_{r+\ell}$ such that $S_{a_j}$ is not $\Gcyc{a_j}$-embeddable in $G_j$ for any $j \in \{1, 2, \ldots, r\}$, and $\Mnest{b_j}$ is not $\Gcyc{b_j}$-embeddable in $G_{r+j}$ for any $j \in \{1, 2, \ldots, \ell\}$. Then every vertex of $G_j$ has degree at most $a_j - 2$, so $|E(G_j)| \leqslant \lfloor \frac{(a_j - 2)n}{2} \rfloor$ for each $j \in \{1, 2, \ldots, r\}$. Moreover, as noted after Lemma~\ref{lemma:nested_reflection}, we have $|E(G_{r+j})| \leqslant \frac{(b_j - 2)n}{2}$ for each $j \in \{1, 2, \ldots, \ell\}$. Therefore,
\begin{equation}\label{eq:star_nested_count}
    \binom{n}{2} = \sum_{j=1}^{r + \ell} |E(G_j)| \leqslant \sum_{j=1}^r \left\lfloor \frac{(a_j - 2)n}{2} \right\rfloor + \sum_{j=1}^\ell \frac{(b_j - 2)n}{2}.
\end{equation}

If $t_1 = 0$, or $t_1$ is odd and $t_2 = 0$, then $n = \sum_{j=1}^r (a_j - 2) + \sum_{j=1}^\ell (b_j - 2) + 2$, so the right-hand side of \eqref{eq:star_nested_count} is at most $\frac{n(n-2)}{2} < \binom{n}{2}$, which is a contradiction. If $t_1$ is even and positive, then $n = \sum_{j=1}^r (a_j - 2) + \sum_{j=1}^\ell (b_j - 2) + 1$ is odd. Thus, for every $j \in \{1, 2, \ldots, r\}$ such that $a_j$ is odd, we have $\lfloor \frac{(a_j - 2)n}{2} \rfloor = \frac{(a_j - 2)n - 1}{2}$, so the right-hand side of \eqref{eq:star_nested_count} equals $\frac{n(n-1)}{2} - \frac{t_1}{2} < \binom{n}{2}$, which is again a contradiction.

It remains to consider the case when $t_1$ is odd and $t_2 = 1$. Then $n = \sum_{j=1}^r (a_j - 2) + \sum_{j=1}^\ell (b_j - 2) + 1$ is even, so the right-hand side of \eqref{eq:star_nested_count} is at most $\binom{n}{2}$. Hence, equality holds in \eqref{eq:star_nested_count}, so $G_j$ is $(a_j - 2)$-regular for each $j \in \{1, 2, \ldots, r\}$, while Lemma~\ref{lemma:nested_reflection} implies that $|E(G_{r+j}) \cap L_n(s)| = \frac{b_j - 2}{2}$ for each $j \in \{1, 2, \ldots, \ell\}$ and $s \in \{0, 1, 2, \ldots, n-1\}$. Let $q \coloneqq \frac{1}{2} \sum_{j=1}^\ell (b_j - 2)$, and observe that $q$ is odd, since $t_2 = 1$. Also, let $Y$ be the spanning subgraph of $K_n$ with the edge set $\bigcup_{j=1}^\ell E(G_{r+j})$. Then every vertex of $Y$ has degree $n - 1 - \sum_{j=1}^r (a_j - 2) = 2q$, and $|E(Y) \cap L_n(s)| = q$ for each $s \in \{0, 1, 2, \ldots, n-1\}$. By computing the sum of $x + y$ over all the edges $\{x, y\} \in E(Y)$ in two ways, we obtain
\begin{equation}\label{eq:endpoint_sum}
    2q \sum_{v=0}^{n-1} v = \sum_{\{x, y\} \in E(Y)} (x + y) \equiv q \sum_{s=0}^{n-1} s \pmod{n}.
\end{equation}
The left-hand side of \eqref{eq:endpoint_sum} equals $q n(n-1) \equiv 0 \pmod{n}$. On the other hand, since $n$ is even, we have $\frac{n(n-1)}{2} = \frac{n}{2}(n-1) \equiv \frac{n}{2} \pmod{n}$, so the right-hand side of \eqref{eq:endpoint_sum} is congruent to $q \cdot \frac{n}{2} \equiv \frac{n}{2} \pmod{n}$, since $q$ is odd. This is a contradiction. Therefore, $\Rcyc(S_{a_1}, \ldots, S_{a_r}, \Mnest{b_1}, \ldots, \linebreak \Mnest{b_\ell}) \leqslant n$.
\end{proof}

For $r = \ell = 1$, Theorem~\ref{thm:star_nested} gives
\[
    \Rcyc(S_a, \Mnest{b}) = \Rdih(S_a, \Mnest{b}) = a + b - 2 - \begin{cases}
        1, & \text{if } a \text{ is odd and } b \text{ is divisible by four},\\
        0, & \text{otherwise},
    \end{cases}
\]
for any $a \geqslant 2$ and any even $b \geqslant 2$. This formula was conjectured in \cite{BaDamSteSto2026} and first proved by Kucheriya et~al.\ \cite[Theorem~3.3]{KuPetrYan2026}. In their proof of the upper bound, they also count the edges of each color, and when $a$ is odd and $b$ is divisible by four, they sum the endpoints of the edges of one color in two ways. For $r = 0$, Theorem~\ref{thm:star_nested} determines the cyclic and dihedral Ramsey numbers of any number of nested matchings, at most one of which has order divisible by four, which was previously known only for two nested matchings \cite[Theorem~4.29]{BaDamSteSto2026}. Finally, for $\ell = 0$, Theorem~\ref{thm:star_nested} is a special case of Theorem~\ref{thm:main}.

\section{Conclusion}\label{sc:conclusion}

We have shown that the cyclic and dihedral Ramsey numbers of separated caterpillars, for any number of colors, coincide with the Ramsey numbers of stars of the same orders, which were determined by Burr and Roberts \cite{BurrRo1973}. We have also proved that Chvátal's formula \cite{Chvatal1977} for the Ramsey numbers of complete graphs versus trees remains valid for the cyclic and dihedral Ramsey numbers of any graph that contains the monotone cycle of the same order versus a separated caterpillar. Moreover, we have determined the cyclic Ramsey numbers of any combination of stars and nested matchings in which at most one nested matching has order divisible by four.

As special cases, our results include several conjectures from \cite{BaDamSteSto2026, DamDjor2026}, all of which were first confirmed by Kucheriya et al.~\cite{KuPetrYan2026}, partially by the same methods. They also determined $\Rcyc(\Mnest{a}, \Mnest{b})$ for all even $a, b \geqslant 2$, which includes the case when $a$ and $b$ are both divisible by four, where $\Rcyc(\Mnest{a}, \Mnest{b}) = a + b - 3$. However, the multicolor case when at least two of the nested matchings have order divisible by four remains open in general, which leads to the following two problems.

\begin{problem}
Determine the cyclic Ramsey numbers $\Rcyc(\Mnest{b_1}, \Mnest{b_2}, \ldots, \Mnest{b_\ell})$ for any even $b_1, b_2, \ldots, b_\ell \geqslant 2$.
\end{problem}

\begin{problem}\label{problem:star_nested}
Determine the cyclic Ramsey numbers $\Rcyc(S_{a_1}, \ldots, S_{a_r}, \Mnest{b_1}, \ldots, \Mnest{b_\ell})$ for any $a_1, a_2, \ldots, a_r \geqslant 2$ and any even $b_1, b_2, \ldots, b_\ell \geqslant 2$.
\end{problem}

Since Theorem~\ref{thm:main} allows each star to be replaced with an arbitrary separated caterpillar of the same order, it is natural to ask what happens when the stars in Problem~\ref{problem:star_nested} are replaced with arbitrary separated caterpillars.

\begin{problem}
Determine the cyclic and dihedral Ramsey numbers $\Rcyc(H_1, \ldots, H_r, \Mnest{b_1}, \ldots, \linebreak \Mnest{b_\ell})$ and $\Rdih(H_1, \ldots, H_r, \Mnest{b_1}, \ldots, \Mnest{b_\ell})$ for any separated caterpillars $H_1, H_2, \ldots, H_r$ and any even $b_1, b_2, \ldots, b_\ell \geqslant 2$.
\end{problem}

\section*{Declarations}

\noindent\textbf{Generative AI and AI-assisted technologies.}
The authors used ChatGPT and Claude for mathematical exploration, including the development of proof strategies and preliminary arguments, and to assist with the writing and editing of this paper. The mathematical arguments presented in this paper were critically reviewed and verified by the authors, who take full responsibility for its content.
\smallskip

\smallskip
\noindent\textbf{Funding.}
I.~Damnjanović was supported by the Ministry of Science, Technological Development and Innovation of the Republic of Serbia, grant number 451-03-34/2026-03/200102, and the Science Fund of the Republic of Serbia, grant \#6767, Lazy walk counts and spectral radius of threshold graphs --- LZWK.
\smallskip

\smallskip
\noindent\textbf{Conflict of interest.}
The authors declare that they have no conflict of interest.
\smallskip

\smallskip
\noindent\textbf{Data and code availability.}
The results are established by mathematical proofs and do not depend on computational datasets or software.

\end{document}